\documentclass[11pt,a4paper]{article}
\usepackage{amsmath,amssymb,amsthm}

\usepackage{hyperref}

\usepackage[dvips]{color}
\definecolor{HAV}{rgb}{0,0,0.5}

\newtheorem{Theorem}{Theorem}[section]
\newtheorem{Proposition}[Theorem]{Proposition}
\newtheorem{Claim}[Theorem]{Claim}
\newtheorem{Corollary}[Theorem]{Corollary}
\newtheorem{Lemma}[Theorem]{Lemma}

\theoremstyle{definition}
\newtheorem{Remark}[Theorem]{Remark}

\newcommand{\begpf}{\noindent{\bf Proof.}\enspace}
\newcommand{\epf}{{\ifhmode\unskip\nobreak\hfil\penalty50 \hskip1em \else\nobreak\fi \nobreak\mbox{}\hfil\mbox{$\square$} \parfillskip=0pt \finalhyphendemerits=0 \par\vskip5pt}}

\newcommand{\Z}{{\mathbb{Z}}}

\begin{document}

\title{Supercongruences for the Catalan-Larcombe-French numbers}

\author{Frazer Jarvis\\{\sl Department of Pure Mathematics}\\{\sl University of Sheffield, Sheffield S3 7RH, U.K.}\\{\tt a.f.jarvis@shef.ac.uk}\\ ~\\Helena A.~Verrill\\{\sl Department of Mathematics}\\{\sl Louisiana State University, Baton Rouge, Louisiana 70803-4918, USA}\\{\tt verrill@math.lsu.edu}}
\date{\today}
\maketitle

\begin{abstract}
In this short note, we develop the Stienstra-Beukers theory of supercongruences
in the setting of the Catalan-Larcombe-French sequence. We also give some
applications to other sequences.
\end{abstract}


\section{Introduction}
The Catalan-Larcombe-French numbers $P_n$
were first defined by Catalan in \cite[Section 9, p. 195]{catalan},
in terms of the ``Segner numbers''.  Catalan stated that these numbers
could be defined by the a recurrence relation:
\begin{equation}\label{eq3}
n^2P_n-8(3n^2-3n+1)P_{n-1}+128(n-1)^2P_{n-2}=0
\end{equation}
for $n\ge2$, with the initial values given by $P_0=1$, $P_1=8$.

Larcombe and French \cite{lf2000} give a detailed account of properties
of the $P_n$, and
obtained \cite[Equations (23) and (35)]{lf2000}
the following formulas for these numbers: 
\begin{equation}\label{eq1}
P_n=2^n\sum_{i=0}^{\lfloor n/2\rfloor}(-4)^i{2(n-i)\choose n-i}^2{n-i\choose i},
\end{equation}
for $n\ge0$.
\begin{equation}\label{eq2}
P_n=\frac{1}{n!}\sum_{p+q=n}{2p\choose p}{2q\choose q}\frac{(2p)!(2q)!}{p!q!}=\sum_{p+q=n}{2p\choose p}^2{2q\choose q}^2\left/{n\choose p}\right..
\end{equation}
These numbers occur in the theory of elliptic integrals
\cite{lf2000},
and there are relations to
the arithmetic-geometric-mean \cite{jlf2003}.
The first few $P_n$ are  $1, 8, 80, 896, 10816, 137728$.
This is sequence 
\htmladdnormallink{A053175}{http://www.research.att.com/~njas/sequences/A053175} in Sloane's database \cite{sloane}.

In an earlier paper~\cite[Theorem 7, p.16]{JLF}, the first author proved the following:
\begin{Proposition}\label{AFJ}
If we write $n=a_da_{d-1}\ldots a_0$ in base $p$, then
$$P_n\equiv P_{a_d}P_{a_{d-1}}\ldots P_{a_0}\pmod{p}.$$
\end{Proposition}

This implies, 
for example, that no $P_n$ is ever divisible by 3 (as none of
$P_0$, $P_1$ or $P_2$ are), or that $P_n$ is divisible by 5 if and only if 
$n$ has a 2 in its base 5 representation. Surprisingly, the stronger result
that the 5-adic valuation $v_5(P_n)$ (i.e., the power of 5 dividing $P_n$) 
is equal to the number of 2s in this 
base 5 representation also seems to be true, but we have no explanation for 
this.

Further, it was observed 
empirically in \cite[Conjectures 3 and 4, p.19]{JLF}
that
\begin{Claim}
\label{claim1a}
Suppose that $p$ is an odd prime, and that $0\le n\le p-1$. Then
\begin{enumerate}
\item $p|P_n$ if and only if $p|P_{p-1-n}$ is divisible by $p$.
\item $p|P_{\frac{p-1}{2}}$ if and only if $p\equiv5\!\!\pmod{8}$ or 
$p\equiv7\!\!\pmod{8}$.
\end{enumerate}
\end{Claim}
In this article, we prove Claim~\ref{claim1a} 
(Corollaries~\ref{cor:claimA} and \ref{cor:claimB} 
and Remark~\ref{rem:converse}) and
furthermore, we prove the following
\begin{Theorem}
\label{thm:supercongruences}
$P_{mp^r}\equiv P_{mp^{r-1}}\pmod{p^r}$.
\end{Theorem}
In fact, the proofs of the results in the claim are mostly entirely elementary,
and follow from the symmetry of the recurrence relation modulo $p$ when $n$ 
is replaced with $p-1-n$, as we explain briefly
in Section~\ref{sec:symofrecrel}  below. 
In Section~\ref{sec:supercong} we use the theory 
developed by Stienstra-Beukers~\cite{SB} and others to prove 
the mod $p^r$ congruences.
In Section~\ref{Granvillesmethod}, we show how to also obtain these congruences
using results of Granville.

\section{Symmetry of the recurrence relation}
\label{sec:symofrecrel}
\subsection{Symmetries for the Catalan-Larcombe-French sequence}
As stated in the Introduction,
the $P_n$ may be defined by the recurrence relation 
$$n^2P_n-8(3n^2-3n+1)P_{n-1}+128(n-1)^2P_{n-2}=0$$
for $n\ge2$, with the initial values given by $P_0=1$, $P_1=8$.
In this section
we will regard this as defining a sequence $P_0,P_1,\ldots,P_{p-1}$ modulo $p$.
(Of course, we cannot determine 
$P_p$ modulo $p$ from this relation owing to the
coefficient of $p^2$ when we put $n=p$. However, 
Proposition~\ref{AFJ} already tells us
that $P_p\equiv P_1\mbox{ (mod $p$)}$, as $p$ is written $10$ in base $p$.)

First shift the variable:
$$(n+1)^2P_{n+1}-8(3n^2+3n+1)P_n+128n^2P_{n-1}=0$$
This recurrence relation has a lot of symmetry when regarded modulo a prime 
number $p$. Indeed, let us write $m=p-1-n$, and reduce modulo $p$. Then
$$m^2P_{p-m}-8(3m^2+3m+1)P_{p-1-m}+128(m+1)^2P_{p-2-m}\equiv0\mbox{ (mod $p$)}.$$
Multiply throughout by $128^m$:
$$128^mm^2P_{p-m}-128^m8(3m^2+3m+1)P_{p-1-m}+128^m128(m+1)^2P_{p-2-m}\equiv0\mbox{ (mod $p$)},$$
and put $Q_m=128^mP_{p-1-m}$. Then
$$128m^2Q_{m-1}-8(3m^2+3m+1)Q_m+(m+1)^2Q_{m+1}\equiv0\mbox{ (mod $p$)}.$$
We conclude that $(Q_n)$ satisfies the same recurrence relation as $(P_n)$
(at least, modulo $p$).

\begin{Lemma}
\label{lem:P_p-1}
$P_{p-1}\equiv(-1)^{\frac{p-1}{2}}\!\!\pmod{p}$ and
$16P_{p-2}\equiv(-1)^{\frac{p-1}{2}}\!\!\pmod{p}$.
\end{Lemma}
\begpf
In expression (\ref{eq2})
$$P_n=\frac{1}{n!}\sum_{r+s=n}{2r\choose r}{2s\choose s}\frac{(2r)!(2s)!}{r!s!},$$
put $n=p-1$; all the terms except that with $r=s=\frac{p-1}{2}$ are clearly 
divisible by $p$. So
\begin{eqnarray*}
P_{p-1}&\equiv&\frac{1}{(p-1)!}{p-1\choose\frac{p-1}{2}}^2\frac{[(p-1)!]^2}{[(\frac{p-1}{2})!]^2}\\
&\equiv&\frac{[(p-1)!]^3}{[(\frac{p-1}{2})!]^6}\mbox{ (mod $p$)}.
\end{eqnarray*}
By Wilson's Theorem, $(p-1)!\equiv-1\mbox{ (mod $p$)}$, and it is easy to see 
that $[(\frac{p-1}{2})!]^2\equiv-(-1)^{\frac{p-1}{2}}$. This gives the result 
for $P_{p-1}$.

To get the result for $P_{p-2}$, we use the recurrence relation
$$(n+1)^2P_{n+1}-8(3n^2+3n+1)P_n+128n^2P_{n-1}=0;$$
put $n=p-1$ and reduce mod $p$:
$$-8P_{p-1}+128P_{p-2}\equiv0\mbox{ (mod $p$)},$$
which gives the value in the statement.
\epf

\begin{Corollary}
\label{cor:claimA}
\begin{equation}
\label{eqn:symmetryinClaimA}
128^nP_{p-1-n}\equiv(-1)^{\frac{p-1}{2}}P_n\pmod{p}
\end{equation}
for all $n$
such that $0\le n\le p-1$. In particular, $p|P_n$ if and only if
$p|P_{p-1-n}$.
\end{Corollary}
\begpf
By definition of the $Q_n$ and Lemma~\ref{lem:P_p-1}
it follows that $Q_0=P_{p-1}\equiv(-1)^{\frac{p-1}{2}}\mbox{ (mod $p$)}$,
and $Q_1=128P_{p-2}\equiv(-1)^{\frac{p-1}{2}}8\mbox{ (mod $p$)}$.
Consequently, mod $p$, the values of $Q_0$ and $Q_1$ are identical to those
of $P_0$ and $P_1$, except for the scaling factor of $(-1)^{\frac{p-1}{2}}$.
Since $(Q_n)$ and $(P_n)$ also satisfy the same recurrence relation, we
conclude that $Q_n\equiv(-1)^{\frac{p-1}{2}}P_n\mbox{ (mod $p$)}$ for all $n$
such that $0\le n\le p-1$. Substituting in the definition of $Q_n$, we deduce
(\ref{eqn:symmetryinClaimA}).
\epf

%

\begin{Corollary}
\label{cor:claimB}
If $p\equiv5\!\!\pmod{8}$ or
$p\equiv7\!\!\pmod{8}$, then $p|P_{\frac{p-1}{2}}$.
\end{Corollary}
\begpf
The central point of the symmetry 
(\ref{eqn:symmetryinClaimA}) 
is when $n=\frac{p-1}{2}$. In this case,
Corollary~\ref{cor:claimA} gives:
$$128^{\frac{p-1}{2}}P_{\frac{p-1}{2}}\equiv(-1)^{\frac{p-1}{2}}P_{\frac{p-1}{2}}\mbox{ (mod $p$)}.$$
So if
 $P_{\frac{p-1}{2}}\not\equiv0\mbox{ (mod $p$)}$
then $(-128)^{\frac{p-1}{2}}\equiv1\mbox{ (mod $p$)}$. This occurs when
$(\frac{-128}{p})=(\frac{-2}{p})=1$, so $-2$ is a quadratic residue, which
means that $p\equiv1\mbox{ (mod 8)}$ or $p\equiv3\mbox{ (mod 8)}$.
The contrapositive gives the result.
\epf

\begin{Remark}
\label{rem:converse}
Using much more sophisticated techniques,
Beukers and
Stienstra \cite{SB} prove that if 
$p\equiv1\mbox{ (mod 8)}$ or $p\equiv3\mbox{ (mod 8)}$, so that $p=a^2+2b^2$ 
for some $a$ and $b$, then
$$P_{\frac{p-1}{2}}\equiv(-1)^{\frac{p-1}{2}}4a^2\mbox{ (mod $p$)}.$$
In particular, $P_{\frac{p-1}{2}}\not\equiv0\mbox{ (mod $p$)}$.
Thus the converse to Corollary~\ref{cor:claimB} also holds.
\end{Remark}

\subsection{Symmetries for other sequences}
Recall that in~\cite{JLF}, we noted that $P_{\frac{p-1}{2}}$ is divisible by $p$
if and only if the {\it Franel number} $f_{\frac{p-1}{2}}$ is divisible by $p$
(indeed, they are congruent modulo~$p$).
Here, $f_n=\sum_{r=0}^n{n\choose r}^3$. 

\begin{Corollary}
Suppose that $p$ is an odd prime. Then
$$\sum_{r=0}^{\frac{p-1}{2}}{\frac{p-1}{2}\choose r}^3\equiv0\!\pmod{p}~\Longleftrightarrow~p\equiv5\!\!\!\pmod{8}\mbox{ or }p\equiv7\!\!\!\pmod{8}.$$
\end{Corollary}

Of course, the method of proof of
Corollary~\ref{cor:claimA}
also applies to other similar recurrence relations, and
it turns out that the Franel numbers furnish another example. As was proven by 
Cusick~\cite{Cu}, the Franel numbers also satisfy a recurrence relation:
$$(n+1)^2f_{n+1}=(7n^2+7n+2)f_n+8n^2f_{n-1},$$
with $f_0=1$, $f_1=2$. In exactly the same way as above, one can prove
\begin{Lemma}
$f_n\equiv(-8)^nf_{p-1-n}\!\!\pmod{p}$ for all $n$
such that $0\le n\le p-1$. In particular, $p|f_n$ if and only if
$p|f_{p-1-n}$.
\end{Lemma}
Indeed, this follows from a similar symmetry argument, once we verify this for
$n=0$ and $n=1$:
$$f_{p-1}=\left({p-1\choose0}^3+{p-1\choose1}^3\right)+\cdots+{p-1\choose p-1}^3$$
and any two consecutive terms in the sum is the sum of two cubes, and is
therefore divisible by the sum of the two numbers. In particular,
$${p-1\choose r}^3+{p-1\choose r+1}^3$$
is divisible by
$${p-1\choose r}+{p-1\choose r+1}={p\choose r+1},$$
which is divisible by $p$. There are an odd number of terms in the sum; pairing
terms off just leaves one term~mod~$p$, $f_{p-1}\equiv{p-1\choose p-1}^3$, say,
so $f_{p-1}\equiv1\equiv f_0\mbox{ (mod $p$)}$. Reducing the recurrence relation
mod~$p$, we get $f_{p-1}+4f_{p-2}\equiv0\mbox{ (mod $p$)}$, so that
$(-8)f_{p-2}\equiv2\equiv f_1\mbox{ (mod $p$)}$. The result then follows with
a similar argument to Corollary~\ref{cor:claimA}.

Other series in Table 7 of
Stienstra and Beukers~\cite{SB} can also be treated in the same
way. We summarise the results:

\begin{Lemma}
\begin{enumerate}
\item Define $a_n=\sum_{k=0}^n{n\choose k}^2{n+k\choose k}$ so that
$$(n+1)^2a_{n+1}=(11n^2+11n+3)a_n+n^2a_{n-1}.$$
Then for any prime $p$, and $0\le n\le p-1$,
$$a_n\equiv(-1)^na_{p-1-n}\pmod{p}.$$
\item Define $b_n=\sum_{k=0}^n{n\choose k}^2{2k\choose k}$ so that
$$(n+1)^2b_{n+1}=(10n^2+10n+3)b_n-9n^2b_{n-1}.$$
Then for any prime $p>3$, and $0\le n\le p-1$,
$$b_n\equiv\left(\frac{-3}{p}\right)9^nb_{p-1-n}\pmod{p}.$$
\end{enumerate}
\end{Lemma}

\section{The Picard-Fuchs equation}
To say more about the $P_n$, we will apply the theory of Beukers and others;
we wish to view the numbers $P_n$ as the coefficients
for a generating function which satisfies a certain differential equation. We
then want to interpret this differential equation as a Picard-Fuchs equation 
for a pencil of elliptic curves.

\begin{Lemma}\label{L1}
The function $\mathbf{P}(x)=\sum_{n=0}^\infty P_nx^n$ is a solution to the
second order differential equation
\begin{equation}
\label{eqn:diffeqnforP_ngenfun}
(1-16x)(1-8x)x\frac{d^2y}{dx^2}+(384x^2-48x+1)\frac{dy}{dx}-8(1-16x)y=0.
\end{equation}
\end{Lemma}
\begpf
This follows easily from the recurrence relation for the numbers $P_n$.
\epf

\begin{Remark}
An alternative proof 
of Lemma~\ref{L1} can be obtained by writing the generating function
$\mathbf{P}(x)$
in terms of a certain elliptic integral
$K(c)$ given in \cite{BB},~\S1.5, and using the 
differential equation for $K(c)$ given in \cite{BB}.
\end{Remark}

We now wish to view (\ref{eqn:diffeqnforP_ngenfun})
as a Picard-Fuchs equation for a
pencil of elliptic curves. 
It turns out that a very similar equation has already appeared in a
paper of the second author \cite{V1}. Indeed, on line~6, Table~6 of \cite{V1}, 
we see that the equation
\begin{equation}
\label{eqn:Picard-Fuchs_eqn}
t(4t-1)(8t-1)f''+(96t^2-24t+1)f'+4(8t-1)f=0\
\end{equation}
is the Picard-Fuchs differential equation for the family of elliptic curves
with level 8 structure, with choice of uniformizing parameter
\begin{equation}
\label{eqn:defoft}
t(\tau)=\frac{\eta(\tau)^4\eta(4\tau)^2\eta(8\tau)^4}{\eta(2\tau)^{10}},
\end{equation}
a weight 0 modular function for $\Gamma_0(8)$,
where $\eta$ is the Dedekind eta function, defined by
$$\eta(\tau)=q^{1/24}\prod_{n\ge 1}(1-q^n),$$
and $q=\exp(2\pi i\tau)$.

The Picard-Fuchs equation is
the equation satisfied by the period of this family of curves, and this
is given by the weight~1 modular form of $\Gamma_0(4)$, given by
\begin{equation}
\label{eqn:defoff}
f(\tau)=\frac{\eta(2\tau)^{10}}{\eta(\tau)^4\eta(4\tau)^4}.
\end{equation}

In view of the above discussion, we have the following
\begin{Theorem}
\label{thm:t-expansion-for-f}
Let $f$ and $t$ be defined by
\begin{eqnarray*}
t(\tau)&=&\frac{1}{2}\frac{\eta(\tau)^4
\eta(4\tau)^2\eta(8\tau)^4}{\eta(2\tau)^{10}}\\
&=&\tfrac{1}{2}q-2q^2+6q^3-16q^4+39q^5-88q^6+188q^7-384q^8+\tfrac{1509}{2}q^9-1436q^{10}+\cdots
\\
f(\tau)&=&\frac{\eta(2\tau)^{10}}{\eta(\tau)^4\eta(4\tau)^4}\\
&=&1+4q+4q^2+4q^4+8q^5+4q^8+4q^9+8q^{10}+8q^{13}+\cdots.
\end{eqnarray*}
Then in a neighbourhood of $\tau=i\infty$, we have
$$\sum_{n\ge 0}P_nt(\tau)^n=f(\tau).$$
\end{Theorem}
\begpf
This follows from the results given in \cite{V1}, and from the fact that
equation 
(\ref{eqn:diffeqnforP_ngenfun})
given in Lemma~\ref{L1} can be obtained from
(\ref{eqn:Picard-Fuchs_eqn})
by setting $y=f$ and $t=2x$.
\epf

\section{Supercongruences via the method of Stienstra-Beukers}
\label{sec:supercong}
Work on the Picard-Fuchs equation by Stienstra and Beukers
\cite{SB}
led to higher
congruences (``supercongruences'') for various quantities defined by similar
recurrence relations. We wish to explore whether there are similar 
supercongruences for the $P_n$ from this general theory.


For convenience, we first give a simple result:
\begin{Lemma}
\label{lem:modppower}
If $t(u)$ is a polynomial in $\Z[u]$, then
for a prime $p$ and integer $k\ge 0$ we have
$$t^{p^k}(u^p)\equiv t^{p^{k+1}}(u)\pmod{p^{k+1}}$$
\end{Lemma}
\begpf
We prove the result by induction on $k$.
\begin{equation}
t(u^p)\equiv t(u)^p\pmod{p}
\end{equation}
i.e., the result holds for $k=0$.
Now suppose 
$t^{p^{k-1}}(u^p)\equiv t^{p^{k}}(u)\pmod{p^k}$
for some $k\ge 1$, i.e.,
$$t^{p^{k-1}}(u^p)= t^{p^{k}}(u) + p^kf(u)$$ 
for some polynomial $f(u)\in\Z[u]$.
Taking $p$th powers of both sides we get
$$t^{p^{k}}(u^p)= (t^{p^{k}}(u) + p^kf(u))^p
=t^{p^{k+1}}(u) + \sum_{i=1}^p \binom{p}{i}p^{ik}t^{p^k(p-i)}(u)f^i(u).
$$ 
When $i=1$, the summand is divisible by $\binom{p}{1}p^k=p^{k+1}$,
and for $i>1$, the summand is divisible by $p^{ik}$, with $ik\ge 2k\ge k+1$,
since $k\ge 1$.  Hence the result follows.
\epf

The next result, following the method of Beukers, is a variant of 
\cite[Proposition~3]{beukersAperynumbers}.

\begin{Proposition}
\label{prop:cong_implies_cong}
Let $t$ be the power series
$$t=\frac{1}{m}\sum_{n\ge 1} a_nu^{n/v},$$
convergent in a neighbourhood of $u=0$,
with $m,v$ positive integers, $a_n\in\Z$ and $a_1=1$.
Suppose that in some neighbourhood of $u=0$ we have an
equality of convergent power series given by
\begin{equation}
\label{eqn:u_and_t_series}
\sum_{n\ge 1}b_n t^{n-1}dt = \sum_{n\ge 1}c_nu^{n-1}du,
\end{equation}
for some integers $b_n$ and $c_n$, $n\ge 1$.

Assume $p$ is a prime not dividing $m$ or $v$.
Then if 
\begin{equation}
\label{eqn:cong1}
c_{mp^r}\equiv c_{mp^{r-1}}\pmod{p^r},
\end{equation}
then we also have
\begin{equation}
\label{eqn:cong2}
b_{mp^r}\equiv b_{mp^{r-1}}\pmod{p^r}.
\end{equation}
\end{Proposition}
\begpf
By \cite[Proposition 1.1]{V2}
the congruence (\ref{eqn:cong1}) is equivalent to 
\begin{equation}
\label{Omegat}
\Omega(t)-\frac{1}{p}\Omega(t^p)=d\theta(t)
\end{equation}
where 
$\Omega(t)=\sum_{n\ge 1}b_n t^{n-1}dt$
and $\theta(t)\in\Z_p[[t]]$,
and (\ref{eqn:cong2}) is equivalent to 
\begin{equation}
\label{Omegau}
\tilde\Omega(u)-\frac{1}{p}\tilde\Omega(u^p)=d\tilde\theta(u)
\end{equation}
where
$\tilde\Omega(u)=\sum_{n\ge 1}c_nu^{n-1}du$
and
$\tilde\theta(u)\in\Z_p[[u]]$.

We can write $\Omega$ and $\tilde\Omega$ as
\begin{eqnarray}
\Omega(t)=df(t) 
&\>\>\text{ and }\>\>&
\tilde\Omega(u)=d\tilde f(u)
\end{eqnarray}
where   $f(t)=\sum_{n\ge 1}\frac{b_n}{n}t^n$ and
$\tilde f(u)=\sum_{n\ge 1}\frac{c_n}{n}u^n$.
By hypothesis, we have $\Omega(t(u))=\tilde\Omega(u)$,
i.e., $df(t(u))=d\tilde f(u)$, so $\tilde f(u)=f(t(u))+\text{\it const}$.
Now we have
\begin{eqnarray*}
\Omega(t^p) -  \tilde\Omega(u^p)
&=&
d(f(t^p(u))-\tilde f(u^p)\\
&=&
d(f(t^p(u))-f(t(u^p))
\\
&=&d\left[\sum_{n\ge 1}\frac{b_n}{n}\Big(t^{np}(u)-t^n(u^p)
\Big)
\right].\\
\end{eqnarray*}
Note that Lemma~\ref{lem:modppower} also applies to polynomials 
in $\Z_p[u]$, and for $p\nmid m$, we have $\frac{1}{m}\in\Z_p$.
Taking limits of sequences of polynomials, 
Lemma~\ref{lem:modppower} also applies to 
power series in $u$ (and fractional powers of $u$) and
in particular to our $t(u)$.  
For any positive integer $n$, write $n=mp^k$ where $(n,m)=1$.
Then replacing $t$ with $t^m$ in Lemma~\ref{lem:modppower}, we get
$$t^{np}(u)\equiv t^n(u^p)\pmod{p^{k+1}},$$
i.e., $t^{np}(u)-t^n(u^p)$ is divisible by
$np$ in $\Z_p$.
Thus 
\begin{equation}
\frac{1}{p}\Omega(t^p) -  \frac{1}{p}\tilde\Omega(u^p) = dg(u).
\end{equation}
where $g(u)\in\Z_p[[u]]$. 
Finally, we need to check that
$d\theta(t)=df(u)$ for some $f(u)\in\Z[[u]]$.
This follows from
inverting the expansion for $t$ in terms of $u$, 
which provided $m,v$ are not divisible by
$p$, and since $a_1=1$, 
gives an expansion of the form $u=\sum \alpha_n t^n$ with
$\alpha_n\in\Z_p$.
From this we have $dt^k = d(\sum \alpha_n u^n)^k=d(\sum \beta_n u^n)$
where $k$ and $\beta_n$ are integers.
Hence (\ref{Omegat}) implies (\ref{Omegau}), and
thus (\ref{eqn:cong1}) implies (\ref{eqn:cong2}).
\epf
In the application of this result, we will take $m=v=2$ and $u=q^2$.

\begin{Lemma}
\label{lem:fandE}
Let $t$ and $f$ be as in (\ref{eqn:defoft}) and (\ref{eqn:defoff}).  Then
$$f\frac{q\frac{dt}{dq}}{t}=1-4q^2-4q^4+32q^6-4q^8-104q^{10}+32q^{12}+192q^{14}+\cdots$$
is an Eisenstein series of weight 3 on $\Gamma_0(8)$ (and a non-trivial
character), and furthermore we have
\begin{equation}
\label{eqn:etaproductforfqdt/tdq}
f\frac{q\frac{dt}{dq}}{t}(\tau)=E(2\tau)
\end{equation}
where
\begin{equation}
\label{eqn:defofE}
E(\tau)=\frac{\eta(\tau)^4\eta(2\tau)^6}{\eta(4\tau)^4}
\end{equation}
\end{Lemma}
\begpf
By \cite[Lemma~0.3]{V2},
$\frac{q\frac{dt}{dq}}{t}$ is a holomorphic modular form of weight 2, so that
$f(q)\frac{q\frac{dt}{dq}}{t}$ is a holomorphic modular form of weight 3.
One can check directly that $t(\tau), f(\tau), E(2\tau)$ are modular forms for
$\Gamma_0(8)$ with a certain character,
using the transformation properties of $\eta$, as given for example in
\cite[Theorem 3.4]{apostol}.  We can alternatively refer to 
the eighth case listed in \cite[Table~1, p.4852]{Martin_eta},
to see that $E(\tau)$ is a Hecke eigenform of weight $3$ and level $4$,
the sixth case in the same table to see that 
$f(\tau)$ is a Hecke eigenform of weight $1$ and level $4$,
and to the last entry in \cite[Table~3, line~17]{CN} to see
that $t(\tau)$ is a weight $0$ modular function for $\Gamma_0(8)$.

Thus  $f(q)\frac{q\frac{dt}{dq}}{t}$ and 
$E(2\tau)$ are modular forms of weight $3$ for $\Gamma_0(8)$, with some
character.
Since the space of weight $3$ modular forms for $\Gamma_1(8)$ is 
finite dimensional,
the equality (\ref{eqn:etaproductforfqdt/tdq}) can be
obtained by comparison
of sufficiently many terms of the $q$-expansions, computed using a
computer program such as PARI, for example.
(One could  be more precise; for example,
we can show that these
forms are modular forms for 
$\Gamma_0(8)\cap\Gamma_1(4)$, for which,
using \cite[Theorem 2.25]{shimura},
the space of weight $3$ modular forms has dimension 4.
One can determine a basis of
Eisenstein series, also given in terms of
eta products, vanishing at all but one of each of the four
cusps, and show that one only needs to compare the
coefficients of $1,q,q^2,q^3$ to determine the equality
(\ref{eqn:etaproductforfqdt/tdq}).  
See also the modular forms given in \cite[Table~11]{ASwD}.)
\epf

\begin{Lemma}
\label{lem:congforcn}
Let $c_n$ be a sequence of integers such that $E(\tau)$ has
$q$-expansion
\begin{equation}
\label{eqn:defofcn}
E(\tau)=1 -4\sum_{n\ge 1}c_nq^n.
\end{equation}
Then
\begin{equation}
\label{eqn:congforcn}
c_{mp^r}\equiv c_{mp^{r-1}}\pmod{p^r}
\end{equation}
\end{Lemma}
\begpf
Fine \cite[p. 85, Eq. (32.7)]{fine} tells us that
$$c_n=\sum_{d|n,d \equiv 1\!\!\!\pmod{4}}d^2-\sum_{d|n,d \equiv 3\!\!\!\pmod{4}} d^2.$$
Thus for a prime $p>2$,
$$c_{mp^r} - c_{mp^{r-1}}=
\sum_{d|m, dp^r \equiv 1\!\!\!\pmod{4}} (dp^r)^2 
-\sum_{d|m, dp^r \equiv 3\!\!\!\pmod{4}} (dp^r)^2\equiv 0\pmod{p^r}.$$
See sequences 
\htmladdnormallink{A120030}{http://www.research.att.com/~njas/sequences/A120030}
and
\htmladdnormallink{A002173}{http://www.research.att.com/~njas/sequences/A002173}
in Sloane's database~\cite{sloane} for further references on the $c_n$.
\epf

\noindent{\bf Proof of Theorem~\ref{thm:supercongruences}.\enspace}
By Theorem~\ref{thm:t-expansion-for-f} we have 
$f(\tau)=\mathbf{P}(t(\tau))$,
and by Lemma~\ref{lem:fandE}, 
and using the expression for $E(\tau)$ given in Lemma~\ref{lem:congforcn},
we have
$f\frac{dt}{t}=E(2\tau)\frac{dq}{q}$, so,
setting $u=q^2$, we have
$$\sum P_n t^{n-1}dt=
\mathbf{P}(t) \frac{dt}{t}=f\frac{dt}{t}=\left(1-4\sum_{n\ge1}
c_n u^n\right)\frac{du}{2u}.$$

Now we apply Proposition~\ref{prop:cong_implies_cong} 
with $u=q^2$, $v=2$ and $m=2$.
We take the $c_n$ 
and $b_n$ of (\ref{eqn:u_and_t_series})
taken to be the 
 $-2c_n$ of (\ref{eqn:defofcn}) (for $n\ge 1$),
and the
Catalan-Larcombe-French numbers $P_n$ respectively.

By Lemma~\ref{lem:congforcn}, with the $c_n$ as defined by
(\ref{eqn:defofE}) and (\ref{eqn:defofcn}),
the congruence (\ref{eqn:congforcn}) holds.  This also holds for $-2c_n$.
Thus the congruence for the $P_n$ follows from
Proposition~\ref{prop:cong_implies_cong}.
\epf

\section{Supercongruences via Granville's method}
\label{Granvillesmethod}

In this section we show how the supercongruences we are interested in
can be obtained in an alternative manner.

We begin by establishing some general results for congruences of binomial
coefficients, following work of Granville~\cite{Gr}. The next result is
Theorem~1 of~\cite{Gr}.

\begin{Theorem}[Granville]\label{Granville}
Suppose that $p^q$ is an odd prime power, and $n=m+r$. Write 
$n=n_dp^d+\cdots+n_0$ in base~$p$,and let $N_j$ be the least residue of
$\lfloor\frac{n}{p^j}\rfloor$ modulo $p^q$ for each $j\ge0$; make corresponding
definitions of $m_j$, $M_j$, $r_j$, $R_j$. Let $e_j$ be the number of indices
$i\ge j$ with $n_i<m_i$ (the number of base~$p$ carries beyond the $j$th digit
in adding $m$ and $r$). Then
\begin{equation}\label{Graneq}
\frac{1}{p^{e_0}}{n\choose m}\equiv(-1)^{e_{q-1}}\left(\frac{(N_0!)_p}{(M_0!)_p(R_0!)_p}\right)\left(\frac{(N_1!)_p}{(M_1!)_p(R_1!)_p}\right)\cdots\left(\frac{(N_d!)_p}{(M_d!)_p(R_d!)_p}\right)\pmod{p^q},
\end{equation}
where $(k!)_p$ denotes the product of the integers $\le n$ not divisible 
by $p$.
\end{Theorem}

(Note that $e_0$ is the number of carries in the base~$p$ sum $m+r=n$, so
this confirms the claims made in the course of the proof of
Proposition~\ref{claim1}.)

Recall that Ljunggren proved the following congruence:
$${pn\choose pm}\equiv{n\choose m}\mbox{ (mod $p^3$)}$$
for $p\ge5$ and any integers $n$ and $m$. In fact, Jacobsthal showed that
this congruence holds modulo~$p^q$, the power of $p$
dividing $p^3mn(m-n)$, and that this is usually best possible: see~\cite{Gr}
for more on this. Using Theorem~\ref{Granville}, we can prove the following:
\begin{Corollary}\label{ljunggren}
Notation as in Theorem~\ref{Granville}. Then 
\begin{equation}\label{Grancoreq}
\frac{1}{p^{e_0}}{pn\choose pm}\equiv\left(\frac{((pN_0)!)_p}{((pM_0)!)_p((pR_0)!)_p}\right).\frac{1}{p^{e_0}}{n\choose m}\pmod{p^q}.
\end{equation}
\end{Corollary}
\begpf
This simply follows on observing that there is only one additional term in
the product when (\ref{Graneq}) is applied with $pn$ and $pm$.
\epf

If $n$ and $m$ are divisible by $p^q$, we can deduce further results. Indeed, 
notice that $(p^q!)_p\equiv-1\mbox{ (mod $p^q$)}$, simply by pairing off a
number less than $p^q$ and not divisible by $p$ with its multiplicative inverse
modulo~$p^q$, leaving only $\pm1$ which are self-inverse (as in one of the
proofs of Wilson's Theorem). It follows that 
$((mp^q)!)_p\equiv(-1)^m\mbox{ (mod $p^q$)}$, for much the same reason. As
a corollary to this observation and Corollary~\ref{ljunggren}, we deduce
the following:
\begin{Corollary}\label{ljunggrensp}
Let $p^q$ be an odd prime power. Then if $p^{e_0}$ denotes the power of $p$
dividing the binomial coefficient ${mp^r\choose kp^s}$, and $r\ge s\ge q$,
then
$$\frac{1}{p^{e_0}}{mp^r\choose kp^s}\equiv\frac{1}{p^{e_0}}{mp^{r-1}\choose kp^{s-1}}\pmod{p^q}.$$
\end{Corollary}
\begpf
Indeed, in the previous corollary, we observe that the numerators and
denominators are of the form $((kp^s)!)_p$ for various values of $k$, and
consequently are all $\pm1$ modulo~$p^s$, and therefore modulo~$p^q$. It is
easy to see that the powers of $-1$ cancel.
\epf

We now return to our study of the Catalan-Larcombe-French numbers. 

We recall (see (\ref{eq1})) that 
$$P_n=2^n\sum_{i=0}^{\lfloor n/2\rfloor}(-4)^i{2(n-i)\choose n-i}^2{n-i\choose i},$$
and we write 
$$g(i,n)={2(n-i)\choose n-i}^2{n-i\choose i}=\frac{((2n-2i)!)^2}{((n-i)!)^3i!(n-2i)!},$$
so that $P_n=2^n\sum_{i=0}^{\lfloor n/2\rfloor}(-4)^ig(i,n)$.

Throughout the rest of this section, we suppose that $p$ is an odd prime.
We first prove that if $p\nmid i$, then
$g(i,mp^r)\equiv0\mbox{ (mod $p^r$)}$. That is, we prove the following
proposition:
\begin{Proposition}\label{claim1}
Let $p$ be an odd prime, and let $i$ and $n$ be integers with $p\nmid i$. Then
$$v_p\left({2(n-i)\choose n-i}^2{n-i\choose i}\right)\ge v_p(n).$$
\end{Proposition}
\begpf
We recall from~\cite[Lemma~2]{JLF}, that
$$v_p\left({s\choose t}\right)=\frac{S_p(t)+S_p(s-t)-S_p(s)}{p-1},$$
where $S_p(s)$ denotes the sum of the digits of $s$ written in base~$p$. 
From~\cite[proof of Lemma~1]{JLF} 
$S_p(t)+S_p(s-t)-S_p(s)=(p-1)c(t,s-t)$, where $c(t,s-t)$ denotes the number
of ``carries'' in the base~$p$ sum $t+(s-t)=s$. It follows that
$v_p\left({s\choose t}\right)$ is exactly $c(t,s-t)$.

In the same way,
$$v_p(g(i,n))=\frac{3S_p(n-i)+S_p(i)+S_p(n-2i)-2S_p(2n-2i)}{p-1}$$
can also be written as $c(n-i,n-i)+c(n-i,i,n-2i)$, the total number of carries
in the two sums:
\begin{eqnarray*}
(n-i)+(n-i)&=&2n-2i;\\
(n-i)+i+(n-2i)&=&2n-2i.
\end{eqnarray*}
Suppose that $p^r|n$, but that $p^{r+1}\nmid n$. Then the base~$p$ expansion
of $n$ ends with $r$ digits~$0$. Since $n-2i$ and $2n-2i$ differ by $n$,
the final $r$ base~$p$ digits of $(n-i)+i$ are all $0$, and the sum in base $p$
$(n-i)+i=n$ will require carries in each of the last $r$ positions (as the
final base~$p$ digit of $i$ is non-zero).
\epf

Perhaps a short illustrative example is in order. Suppose that $n=18$, $p=3$
and $i=7$. Then, in base~3:
$$n=200,\quad i=21,\quad n-i=102,\quad n-2i=11,\quad 2n-2i=211.$$
The number of carries in $(n-i)+(n-i)=(2n-2i)$ is the number of carries in
$102+102=211$, which has one carry. More importantly, the other sum is
$102+21+11=211$, and the number of carries is the same as that in the sum
$102+21=200$. However, to end with two zeros will require two carries, and,
in fact, there are exactly two carries. We see that 
$v_3(g(7,18))=1+2=3$.

By the proposition, $p^r|g(i,mp^r)$ if $p\nmid i$. In other words,
$g(i,mp^r)\equiv0\mbox{ (mod $p^r$)}$ if $p\nmid i$.
(Indeed, for future reference, the same argument as in the proof
of Proposition~\ref{claim1} tells us that $p^{r-s}|g(kp^s,mp^r)$.)

We conclude that
\begin{eqnarray*}
P_{mp^r}&=&2^{mp^r}\sum_{i=0}^{\lfloor mp^r/2\rfloor}(-4)^ig(i,mp^r)\\
&\equiv&2^{mp^r}\sum_{{i=0\atop p|i}}^{\lfloor mp^r/2\rfloor}(-4)^ig(i,mp^r)\pmod{p^r}\\
&\equiv&2^{mp^r}\sum_{j=0}^{\lfloor mp^{r-1}/2\rfloor}(-4)^{jp}g(jp,mp^r)\pmod{p^r}\\
&\equiv&2^{mp^{r-1}}\sum_{j=0}^{\lfloor mp^{r-1}/2\rfloor}(-4)^{jp}g(jp,mp^r)\pmod{p^r},
\end{eqnarray*}
the last congruence following from the Fermat-Euler Theorem:
$a^{\phi(p^r)}\equiv1\mbox{ (mod $p^r$)}$ for $p\nmid a$~--~apply this with
$a=2^m$, and recall that $\phi(p^r)=p^r-p^{r-1}$.

The result will follow from a consideration of the terms $g(jp,mp^r)$. To
analyse these terms, we will use the congruence results we established earlier,
using the results of Granville~\cite{Gr}.

We now prove the following result.
\begin{Theorem}
$g(jp,mp^r)\equiv g(j,mp^{r-1})\pmod{p^r}$.
\end{Theorem}
\begpf
Let us write $j=kp^s$, with $p\nmid k$. We have already remarked that
$p^{r-s-1}|g(kp^{s+1},mp^r)$, and also that $p^{r-s-1}|g(kp^s,mp^{r-1})$,
since the number of carries in one of the sums (as in Proposition~\ref{claim1})
is at least $r-s-1$. Let $e_0$ denote the total number of carries, as above,
so $e_0\ge r-s-1$. We now want to see that 
\begin{equation}\label{need1}
\frac{1}{p^{e_0}}g(kp^{s+1},mp^r)\equiv\frac{1}{p^{e_0}}g(kp^s,mp^{r-1})\pmod{p^{s+1}},
\end{equation}
i.e.,
\begin{equation}\label{need2}
\frac{1}{p^{e_0}}{2mp^r-2kp^{s+1}\choose mp^r-kp^{s+1}}^2{mp^r-kp^{s+1}\choose kp^{s+1}}\equiv\frac{1}{p^{e_0}}{2mp^{r-1}-2kp^s\choose mp^{r-1}-kp^s}^2{mp^{r-1}-kp^s\choose kp^s}\pmod{p^{s+1}}.
\end{equation}
However, this now follows as in Corollary~\ref{ljunggrensp}, with $q=s+1$.
\epf

This is not quite what we need, as we need to account for the power of $-4$.
Luckily, this is now also manageable:

\begin{Proposition}
$(-4)^{jp}g(jp,mp^r)\equiv(-4)^jg(j,mp^{r-1})\pmod{p^r}$.
\end{Proposition}
\begpf We start by observing that if $j=kp^s$ with $p\nmid k$, then
\begin{eqnarray*}
g(jp,mp^r)=g(kp^{s+1},mp^r)&\equiv&g(kp^s,mp^{r-1})\pmod{p^r}\\
&\equiv&g(kp^{s-1},mp^{r-2})\pmod{p^{r-1}}\\
&\dots&\\
&\equiv&g(k,mp^{r-s+1})\pmod{p^{r-s-1}}\\
&\equiv&0\pmod{p^{r-s-1}}
\end{eqnarray*}
Consequently, $p^{r-s-1}|g(jp,mp^r)$, and again using the proposition,
$p^{r-s-1}|g(j,mp^{r-1})$. Further, again by the Fermat-Euler
Theorem, we have 
$(-4)^{p(kp^s)}=(-4)^{kp^{s+1}}\equiv(-4)^{kp^s}\mbox{ (mod $p^{s+1}$)}$,
and so $p^{s+1}|(-4)^{jp}-(-4)^j$. Then
\begin{eqnarray*}
&&(-4)^{jp}g(jp,mp^r)-(-4)^jg(j,mp^{r-1})\\
&=&(-4)^{jp}g(jp,mp^r)-(-4)^{jp}g(j,mp^{r-1})+(-4)^{jp}g(j,mp^{r-1})-(-4)^jg(j,mp^{r-1})\\
&=&(-4)^{jp}\left(g(jp,mp^r)-g(j,mp^{r-1})\right)+\left((-4)^{jp}-(-4)^j\right)g(j,mp^{r-1})
\end{eqnarray*}
which is divisible by $p^r$ as required.
\epf

We can now deduce our main theorem, giving supercongruences for the
Catalan-Larcombe-French numbers:

\noindent{\bf Proof of Theorem~\ref{thm:supercongruences}.\enspace}
We have already explained that
$$P_{mp^r}\equiv 2^{mp^{r-1}}\sum_{j=0}^{\lfloor mp^{r-1}/2\rfloor}(-4)^{jp}g(jp,mp^r)\pmod{p^r},$$
and now we know that
$(-4)^{jp}g(jp,mp^r)\equiv(-4)^jg(j,mp^{r-1})\mbox{ (mod $p^r$)}$. It
follows that
$$P_{mp^r}\equiv 2^{mp^{r-1}}\sum_{j=0}^{\lfloor mp^{r-1}/2\rfloor}(-4)^jg(j,mp^{r-1})\pmod{p^r}.$$
But the right-hand side is one of the ways to define $P_{mp^{r-1}}$, and the
result follows.
\epf


\begin{thebibliography}{99}

\bibitem{apostol}
T.M.Apostol, Modular functions and Dirichlet series in number theory. 
Graduate Texts in Mathematics 41, Springer-Verlag, New York-Heidelberg, 
1976. 

\bibitem{beukersAperynumbers}
F.Beukers, 
Another congruence for the Ap\'ery numbers, J. Number
Theory \textbf{25} (1987) 201--210

\bibitem{BB} J.M.Borwein, P.B.Borwein, Pi and the AGM, 
Wiley-Interscience (1987)

\bibitem{catalan}
E.Catalan, Sur les Nombres de Segner, Rend. Circ. Mat. Pal.
\textbf{1} (1887) 190--201

\bibitem{CN}
J.H.Conway, S.P.Norton, Monstrous moonshine.
Bull. London Math. Soc. \textbf{11} (1979) 308--339

\bibitem{Cu} T.W.Cusick, Recurrences for sums of powers of binomial
coefficients, J. Combin. Theory Ser.A \textbf{52} (1989) 77--83

\bibitem{ASwD}
Liqun Fang, J.W.Hoffman, B.Linowitz, A.Rupinski, H.A.Verrill,
Modular forms on noncongruence subgroups and Atkin-Swinnerton-Dyer
relations.  arXiv:math.NT/0805.2144 

\bibitem{fine} N.J.Fine, 
Basic hypergeometric series and applications.
With a foreword by George E. Andrews. Mathematical Surveys and Monographs, 27. American Mathematical Society, Providence, RI, 1988. xvi+124 pp. ISBN: 0-8218-1524-5 

\bibitem{Gr} A.Granville, Arithmetic properties of binomial coefficients I.
Binomial coefficients modulo prime powers. Organic mathematics (Burnaby, BC,
1995), 253--276, CMS Conf. Proc., 20, Amer. Math. Soc., Providence, RI, 1997

\bibitem{jlf2003}
A.F.Jarvis, P.Larcombe, D.French,
Applications of the A.G.M. of Gauss: some new properties of the
Catalan-Larcombe-French sequence. Proceedings of the Thirty-Fourth
Southeastern International Conference on Combinatorics, Graph Theory
and Computing.  
Congr. Numer.  \textbf{161} (2003) 151--162

\bibitem{JLF} A.F.Jarvis, P.Larcombe, D.French, On small prime divisibility of
the Catalan-Larcombe-French numbers, Indian Journal of Mathematics \textbf{47}
(2005) 159--181

\bibitem{lf2000}
P.Larcombe, D.French,
On the 'other' Catalan Numbers: A Historical Formulation Re-examined,
Congr. Numer. \textbf{143} (2000) 33--64

\bibitem{Martin_eta}
Y.Martin,  Multiplicative $\eta$-quotients, 
Trans. Amer. Math. Soc. \textbf{348} (1996) 4825--4856

\bibitem{shimura}
G.Shimura, Introduction to the arithmetic theory of automorphic functions.
Reprint of the 1971 original. Publications of the Mathematical Society of 
Japan, 11. Princeton University Press, Princeton, NJ, 1994. xiv+271 pp. ISBN: 0-691-08092-5

\bibitem{sloane}
N.J.A.Sloane, The On-Line Encyclopedia of Integer Sequences,
published electronically at www.research.att.com/$\sim$njas/sequences/.

\bibitem{SB} J.Stienstra, F.Beukers, On the Picard-Fuchs equation
and the formal Brauer group of certain elliptic $K3$-surfaces,
Math. Ann. \textbf{271} (1985) 271--304

\bibitem{V1} H.A.Verrill, Picard-Fuchs equations of some families of
elliptic curves, Proceedings on Moonshine and related topics 
(Montr\'eal, Qu\'ebec, 1999), 253--268, CRM Proc. Lecture Notes, 30, 
Amer. Math. Soc., Providence, RI, 2001.

\bibitem{V2} H.A.Verrill, Some congruences related to modular forms,
Max Planck Institut f\"ur Mathematik preprint \textbf{26} (1999)
\end{thebibliography}
\end{document}